\documentclass[12pt]{article}
\usepackage{amsmath,amssymb,amsthm,amsfonts,mathtools,mathrsfs}
\usepackage{etaremune, enumerate, float, verbatim}
\usepackage{algorithm}
\usepackage{algorithmic}
\usepackage{hyperref, theoremref}
\usepackage{soul} 
\usepackage{graphicx}
\usepackage{url}
\usepackage[mathlines]{lineno}
\usepackage{dsfont} 
\usepackage{tikz, graphicx, subcaption, caption}
\usepackage[algo2e,ruled,vlined]{algorithm2e}
\usepackage{xcolor}
\numberwithin{figure}{section}
\usepackage{theoremref}
\usepackage{listings}
\usepackage[noadjust]{cite}
\usetikzlibrary{positioning}
\usepackage{comment}
\usepackage{bm}

\newtheorem{thm}{Theorem}[section]
\newtheorem{cor}[thm]{Corollary}
\newtheorem{lem}[thm]{Lemma}
\newtheorem{prop}[thm]{Proposition}

\newtheorem{obs}[thm]{Observation}

\newtheorem{quest}[thm]{Question}
\theoremstyle{definition}

\newtheorem{defn}[thm]{Definition}
\newtheorem{ex}[thm]{Example}

\newcommand{\R}{\mathbb{R}}

\newcommand{\msym}{\operatorname{\mathsf{Sym}}} 
\newcommand{\csym}{\operatorname{\mathsf{CSym}}} 
\newcommand{\fort}{\operatorname{\mathsf{Fort}}} 

\newcommand{\calC}{\mathcal{C}}
\newcommand{\calF}{\mathcal{F}}
\newcommand{\calS}{\mathcal{S}}

\newcommand{\calFmin}{\mathcal{F}^{\rm min}}
\newcommand{\calSmin}{\mathcal{S}^{\rm min}}

\DeclareMathOperator{\rank}{rank}

\DeclareMathOperator{\nul}{nul}
\DeclareMathOperator{\Nul}{Nul}
\DeclareMathOperator{\ft}{ft}

\DeclareMathOperator{\supp}{supp}

\DeclareMathOperator{\Z}{Z}
\DeclareMathOperator{\Y}{T}
\DeclareMathOperator{\M}{M}
\DeclareMathOperator{\N}{N}

\title{Compatible Forts and Maximum Nullity of a Graph}
\date{\today}
\author{Veronika Furst \thanks{Department of Mathematics, Fort Lewis College, Durango, CO, USA. (furst\_v@fortlewis.edu)}
\and
John Hutchens \thanks{Department of Mathematics and Statistics, University of San Francisco, San Francisco, CA, USA. (jhutchens@usfca.edu)} \and Lon Mitchell \thanks{Department of Mathematics \& Statistics, Eastern Michigan University, Ypsilanti, MI, USA.  (lmitch50@emich.edu)} 
\and
Yaqi Zhang\thanks{Fundamental Software Research Center, Peking University Chongqing Research Institute of Big Data, Chongqing, China. (yaqizhangus@outlook.com)}}

\begin{document}
\maketitle

\begin{abstract}
We consider bounds on maximum nullity of a graph via transversal numbers of compatible collections of forts.    Results include generalizations of  theorems from symmetric to combinatorially symmetric matrices, special bases of matrix nullspaces derived from transversal sets,   and examples of issues that arise when considering only minimal forts and how to avoid them.  We also show an important difference between constructing symmetric and combinatorially symmetric matrices associated to a graph whose nullspaces are supported on collections of disjoint forts.
\end{abstract}
\medskip

\medskip
\section{Introduction}

Many characteristics of matrices with symmetric zero/nonzero patterns are determined by their associated graphs.  In what follows, we discuss properties of a graph that provide information about the structure of the nullspaces of matrices associated to the graph, including bounds on the maximum nullity of such matrices.

We consider graphs $G=(V(G),E(G))$ that are simple (no multiple edges and no edges from a vertex to itself), and we will assume a knowledge of basic graph theory definitions (see, for example, ~\cite{west}). 
We consider square matrices $A = [a_{ij}]$ with real entries that are either symmetric ($a_{ij} = a_{ji}$) or combinatorially symmetric ($a_{ij}\neq 0$ if and only if $a_{ji}\neq 0$).  For an $n$-by-$n$ combinatorially symmetric matrix $A$, the graph of $A$, denoted $G(A)$,  is the graph with vertices $\{1,2,\dotsc,n\}$ and edges $\{ij: \  i\neq j \mbox{ and } a_{ij}\neq 0\}$. 
For a graph $G$, let $\csym(G)$ be the set of combinatorially symmetric matrices whose graph is $G$ and $\msym(G)$ be the subset consisting of symmetric matrices whose graph is $G$; we say these matrices are ``associated to'' or ``described by'' the graph.  Further define $\N(G)$ to be the maximum nullity among matrices in $\csym(G)$ and $\M(G)$ the maximum nullity among matrices in $\msym(G)$.  We denote the null space and nullity of a matrix $A$ by $\Nul(A)$ and $\nul(A)$, respectively. 
The \emph{support} of a vector $x = (x_i)$, denoted by $\supp(x)$, is the set of indices $i$ for which $x_i \neq 0$. 
 
{\em Zero forcing} is a color changing process on the vertices of a graph~\cite{HogbenLinShaderbook}.  To start, each vertex is either filled with color or unfilled, and each type of zero forcing follows a specific color change rule which can change the color of an unfilled vertex to filled. The process stops when no more vertices can be filled.  The \emph{standard zero forcing} color change rule is to change the color of an unfilled vertex $w$ to filled if $w$ is the unique unfilled neighbor of a filled vertex $v$; we say $v$ {\em forces} $w$. A {\em standard zero forcing set} of a graph $G$ is a set of vertices $B$ such that if $B$ is the set of initially filled vertices and the standard zero forcing color change rule is applied repeatedly, then all vertices of $G$ eventually become filled.   The \emph{standard zero forcing number} of a graph $G$ is $\Z(G)=\min \{|B|:\ B \text{ is a zero forcing set of }G\}$. Zero forcing was first introduced to provide an upper bound for $\M(G)$ \cite{AIMMINRANK}.
The name refers to the fact that if a null vector $x$ of $A\in\msym(G)$ is such that $\supp(x)$ is disjoint from a zero forcing set of $G$, then $x$ must be the $0$ vector.

A \emph{fort} in a graph $G$ is a nonempty subset $F$ of vertices of $G$ such that if $v$ is a vertex of $G$ not in $F$, then $v$ does not have exactly one neighbor in $F$ \cite{FastHicks18}.  Forts resist zero forcing in that no vertex outside of $F$ can force a vertex inside of $F$ to become filled.   Moreover, it is known that $\Z(G)$ equals the transversal number of the collection of all forts of $G$ (see Proposition \ref{Z-tau} and the subsequent comment).

We are primarily interested in extending recent work of Hicks~et~al.~\cite{louis}, who proved that supports of null vectors of symmetric matrices associated with a graph form special subcollections of forts, called \emph{compatible}, and showed that transversal numbers of these special subcollections give a better upper bound for the maximum nullity of such matrices than the zero forcing number.  

In Section 2 we make connections between forts and matroid theory. In particular we show that any finite matroid corresponds to some collection of forts of a finite simple graph.  This provides a new approach to studying graph parameters using matroids associated to various collections of forts for a given graph.  In Section 3 we address these connections and consider that we can use transversals of compatible collections of forts to identify useful bases of nullspaces of matrices.  We also show that much of the work of Hicks~et~al.~extends directly to combinatorially symmetric matrices.  

In Section 4 we show that while Hicks~et~al.~elected to work with minimal forts, doing so can cause problems for their desired results.  We provide examples of these issues and show how to avoid them.  We introduce a graph parameter $\Y(G)$ such that $\M(G) \leq \Y(G) \leq \Z(G)$ for any graph $G$ and show that there are examples where $\Y(G)< \Z(G)$ and $\M(G) < \Y(G)$.  This is related to Question 2 in \cite{AIMMINRANK}, which seeks a tighter upper bound for $\M(G)$ than $\Z(G)$.  It also establishes a matroid theoretic parameter related to $\M(G),\N(G)$, and $\Z(G)$. 

Although much of our work is equally valid for symmetric and combinatorially symmetric matrices associated to a given graph, our final section discusses how these two matrix classes differ when we try to construct a matrix that has a given subcollection of forts as its null vector supports.  In particular, for subcollections of disjoint forts, we show that it is not always possible to construct a symmetric matrix with those forts as its null vector supports and characterize when it is possible.  

\section{Forts and Matroids}\label{Section:FortsAndMatroids}

  Given a graph $G$, let $\fort(G)$ be the collection of all forts of $G$.  We denote the set of nontrivial supports of null vectors of a matrix $A\in \csym(G)$ by
\[ \calF_A = \{ \supp(x): \ x \in \Nul(A), \ x \neq 0 \}. \]
If $\calS$ is a collection of sets, then we denote the collection of minimal sets in $\calS$ by $\calSmin$, where a set is {\em minimal} if it does not contain any other set in $\calS$.  In particular, the collections of minimal forts in $G$ and minimal null vector supports are denoted by $\fort(G)^{\min}$
and $\calF_A^{\min}$, respectively.  We will often use similar notation for forts and null vector supports due to the connection that for a nonzero vector $x \in \R^{|G|}$, $Ax = 0$ for some $A \in \msym(G)$ if and only if and only if $\supp(x)$ is a fort of $G$:

\begin{thm}[\cite{louis,DFFHMZ}]\label{FortSupport}
    For any matrix $A \in \msym(G)$, the support of any nonzero null vector of $A$ is a fort of $G$. Conversely, for any fort $F$ of $G$ and any vector $x$ whose support is $F$, there is a matrix $A \in \msym(G)$ that has $x$ as a null vector.  
\end{thm}

The proof by Deaett~et~al.~\cite{DFFHMZ} only uses the zero/nonzero pattern of $G$, so we have the following generalization as well.

\begin{prop} \label{csym forts}
    For $0\neq x \in \R^{|G|}$, if $A \in \csym(G)$ and $Ax = 0$, then $\supp(x)$ is a fort of $G$.
\end{prop}

Forts are also connected to matroids via null vector supports.  While there are several equivalent definitions of a matroid, we wish to focus on the sets of vertices in a graph $G$ that correspond to dependent sets of columns of matrices in $\msym(G)$ and $\csym(G)$. This leads us to define a matroid in terms of its minimal dependents sets.  

\begin{defn}[{\cite{Oxley}}] 
\label{def:matroid}
Let $S$ be a finite set, and let $\calC$ be a collection of subsets of $S$.  Then $\bm{M}=(S,\calC)$ is a \emph{matroid} if
\begin{enumerate}
    \item $\varnothing \not\in \calC$; 
    \item if $C_1$, $C_2 \in \calC$ and $C_1 \subseteq C_2$, then $C_1 = C_2$;
    \item if $C_1,C_2 \in \calC$ such that $C_1 \neq C_2$ and $x \in C_1 \cap C_2$, then there is another set $C_3 \in \calC$ such that $C_3 \subseteq (C_1 \cup C_2) \setminus \{x\}$. 
\end{enumerate}
\end{defn} 

An \emph{independent set} of a matroid $\bm{M} = (S,\calC)$ is defined to be any set $I\subseteq S$ that contains no set from $\calC$.  A set $B \subset S$ is a basis for $\bm{M}$ if it is an independent set of maximum size.  The {\em rank} of a matroid $\bm{M}$ is defined to be the size of a basis for $\bm{M}$.  The members of $\calC$ are therefore the minimal {\em dependent} sets of $\bm{M}$, known as the {\em circuits} of $\bm{M}$.  Circuits with one element are called {\em loops}.

For any matrix $A$ there exists a matroid $\bm{M}(A) = (S, \calC)$ where $S$ is the set of columns of $A$ and the minimal sets of dependent columns are the circuits.  If a matroid $\bm{M} = \bm{M}(A)$ for some matrix $A$, we say that $\bm{M}$ is {\em representable}. 
When $\bm{M}$ is representable, $\rank(A) = \rank(\bm{M}(A))$. 

We are interested in (square) matrices that have a zero/nonzero pattern that corresponds to a graph $G$ and thus matrices $A \in \csym(G)$ or $A \in \msym(G)$ such that $(S,\calC) = \bm{M}(A)$.  In either case, by Theorem \ref{FortSupport}, $\calC = \calFmin$ for some $\calF \subseteq \fort(G)$.  For any $A\in \csym(G)$, $\calF_A$ satisfies Properties (1) and (3) of Definition \ref{def:matroid} while $\calF_A^{\min}$ satisfies all three properties.  So we have the following observation.

\begin{obs} \label{F_A^min-matroid}
    If $A \in \csym(G)$, then $\calFmin_A$ comprises the circuits of the matroid $\bm{M}(A) = (V(G), \calFmin_A)$. 
\end{obs} 

Nelson~\cite{AlmostAllMatroids} showed that as $n$ goes to infinity, the proportion of representable matroids approaches zero.  There are classes of matroids known to be representable over certain fields \cite{TransversalRep}. In general the problem of whether or not a matroid is representable over $\mathbb{R}$ is unknown and difficult to determine \cite{MWN}.
On the other hand, the next result shows that any matroid can arise as a matroid whose circuits are collections of forts of a graph $G$.

\begin{prop} \label{AllMatroidForts}
    For every matroid $\bm{M} = (S,\calC)$ there exists a graph $G$ with $\calC \subseteq \fort(G)$ and $S = V(G)$.  In particular, if every circuit in $\calC$ has at least two elements, then $\calC \subseteq \fort(K_n)$ where $n = |S|$.
\end{prop}
\begin{proof}
    Let $\bm{M} = (S,\calC)$ be a matroid with $|S| =n$ and circuits $\calC$.  Assume $\bm{M}$ has $m$ loops (circuits that contain only one element).  If $m\geq 1$, we label the loops in $\calC$ by $\{1,2,\ldots, m\}$.  No other circuit in $\calC$ can contain any loop element since that would contradict the minimality of the sets in $\calC$.  Notice that the complete graph on $d$ vertices satisfies $\fort(K_d) = \{ F \subseteq V(K_d): \  |F| \geq 2 \}$.  All circuits in $\calC$ that are not loops are forts of the complete graph on vertices labeled $\{m+1, m+2, \cdots, n\}$ with $d=n-m$.  So $\calC \subseteq \fort(G)$ for $G = mK_1 \sqcup K_{n-m}$, where we obtain the loops in the matroid from the isolated vertices of $G$.  If $\calC$ contains no loops, then $m=0$ and $\calC \subseteq \fort(K_n)$.
\end{proof}

Finally, forts are also directly connected to zero forcing via minimum transversals: The \emph{transversal number} of a collection of sets $\calS$ is 
\[
\tau(\calS) = \min \{|T|: \ \text{$T \cap S \neq \emptyset$ for all $S \in \calS$}\}.
\]

\begin{prop}[{\cite[Theorem~2.2]{cameron2023forts}}] \label{Z-tau}
For any graph $G$, $\Z(G) = \tau(\fort(G)^{\min})$.
\end{prop}

In fact, the original proof of this result establishes that $\Z(G) = \tau(\fort(G)^{\min}) = \tau(\fort(G))$.  We also note here that in general $\tau(\calS) = \tau(\calSmin)$. 

\section{Compatible Forts and Nullity}\label{Section:CompFortsAndNullity}

Given that the zero forcing number of a graph can be expressed as the transversal number of its collection of (minimal) forts, and also that forts are supports of null vectors, a natural question is whether the maximum nullity of a graph could be related to the transversal number of its forts. Hicks et al.~\cite{louis} were able to connect the nullity of a single matrix to a transversal number.

\begin{prop}[{\cite[Lemma~2]{louis}}]
    For $A \in \msym(G)$, $\nul(A)  =  \tau(\calFmin_A)$.
\end{prop}
The original proof of this result does not rely on either the symmetry of $A$ or the minimality of the supports of its null vectors.  The next lemma states this generalized result; we include its proof for completeness.
\begin{lem} \label{tau(C_A)}
    For $A \in \csym(G)$, $\nul(A)  = \tau(\calFmin_A) = \tau(\calF_A)$. 
\end{lem}
\begin{proof}
    A set $T$ intersects every fort in $\calF_A$ if and only if the vertices in $V(G)\setminus T$ correspond to linearly independent columns of $A$.  The result follows by noting that $|T|$ is minimized precisely when $|V(G)\setminus T|$ is maximized.
\end{proof}
Hicks et al.~\cite{louis} define a collection $\calF$ of forts to be \emph{compatible} if for any distinct $F_1, F_2 \in \calF$ and $x \in F_1 \cap F_2$, there exists $F \in \calF$ such that $F \subseteq (F_1 \cup F_2) \setminus \{x\}$. The definition of compatibility mirrors the exchange property for circuits in a matroid as in Definition \ref{def:matroid}(3).  Indeed, the forts in a compatible collection $\calF \subseteq \fort(G)$ that are minimal with respect to that collection define the circuits of a matroid on the set of vertices.

As a result, compatible collections of forts enjoy other properties similar to those that arise in the study of matroids.  For example, Lemma~\ref{tau(C_A)} can be thought of as a generalization of the following matroid result, which has a similar proof.

\begin{prop}\label{theoneweremovedthenputback}
    Let ${\bm M}=(S,\calC)$ be the matroid defined by $\calC$ as its circuits.  Assume its base set $S$ has cardinality $n$.  Then $\rank(\bm{M}) = n - \tau(\calC)$.
\end{prop}

For another example, if $I$ is an independent set in a matroid but $I\cup\{x\}$ is dependent, then there exists a unique \emph{fundamental circuit} $C(I,x)$ such that $x \in C(I,x) \subseteq I\cup\{x\}$.  A compatible collection of forts exhibits a related property with respect to vertices of a minimum transversal:

\begin{thm}\label{correspondence}
Given a compatible collection of forts \(\calF\) and a minimum transversal set \(T\) of $\calF$, for each vertex $v\in T$ there exists a unique fort \(F_v \in \calF\) such that \(F_v\) contains no other element of \(T\).
\end{thm}

\begin{proof}
There has to be at least one such fort as otherwise \(T\) is not minimum (since we could remove \(v\)).  If there are two such forts, say \(F_1\) and \(F_2\), then by compatibility requirement, there exists \(F\in\calF\) with \(F \subseteq (F_1 \cup F_2) \setminus \{v\}\).  But then \(F\) would not be transversed by \(T\).  
\end{proof}

One example of the  utility of the previous result can be found below in the argument of  Example \ref{ex:petersen}.  As another application, given any combinatorially symmetric matrix, we can find a set of vertices in its graph that correspond to a nice basis of its nullspace:

\begin{prop}\label{specialbasis}
    Let $A\in \csym(G)$, and let $T = \{v_1,v_2,\ldots,v_{n-k}\}$ be a minimum transversal of $\calF_A$.  If $\{F_{v_i}: \ v_i\in T\} \subseteq \calF_A$ is the set of unique forts associated to $T$ as in Theorem \ref{correspondence}, then there exists a basis of $\Nul(A)$ of the form $X=\{x_1,x_2,\ldots,x_{n-k}\}$ such that $\supp(x_i) = F_{v_i}$ for $1 \leq i \leq n-k$.
\end{prop}
\begin{proof}
    By Lemma \ref{tau(C_A)}, $|T| = \tau(\calF_A) = n-k$, where $k= \rank(A)$.  For each $i$, there is a vector $x_i \in \Nul(A)$ with $\supp(x_i) = F_{v_i}$.  Let $X = \{x_i: \ 1\leq i \leq n-k\}.$  By Theorem \ref{correspondence}, for each $i$, $x_i$ is the unique vector in $X$ such that $v_i \in \supp(x_i)$.  So $x_i$ is the unique vector in $X$ with a nonzero $i$th component, and hence $X$ is a linearly independent set.  Since $X\subseteq \Nul(A)$ with $n-k = |X| = \dim(\Nul(A))$, $X$ must be a basis for $\Nul(A)$. 
\end{proof}

Proposition \ref{specialbasis} can prove useful in constructing symmetric matrices of a given nullity and graph. It also has implications for the ``nullspace representations" of the vertices $V(G)$ of the graph $G(A)$ of a symmetric matrix $A$, which arise by selecting a basis $\{v_1,v_2,\dotsc,v_m\}$ for the nullspace of $A$, creating the matrix $B=\left[ \begin{matrix} v_1 & v_2 & \dotsb & v_m\end{matrix}\right]$ with those basis vectors as the columns, and then assigning the rows of $B$ to the vertices of the graph~\cite{nullreps, Lovasz2017}.

\begin{cor}\label{thecor}
For every (combinatorially) symmetric matrix $A$ with graph $G$, there exists a nullspace representation of $V(G)$ corresponding to $A$ that includes the standard basis. 
\end{cor}

In particular, Corollary~\ref{thecor} means that the vertices of $G$ can be (re)labeled so that the top square submatrix of $B$ is an identity matrix.

\section{Problems with Minimality}\label{Section:ProbMinimality}

A minimal fort of $G$, that is, an element of $\fort(G)^{\min}$, is a minimal null support of some $A\in\csym(G)$.  However, a minimal null support of some $A\in\csym(G)$ need not be a minimal fort (since a properly contained fort may be the support of a null vector of a different matrix).  For example, the Laplacian matrix $L$ of a connected graph has a one-dimensional nullspace, spanned by the all-ones vector; therefore, $V(G)$ is a minimal support for $L$, but it clearly need not be a minimal fort.  Hicks et al.~\cite[Theorem~1]{louis} claim that $\M(G) \leq \max \{\tau(\calF):  \calF \subseteq \fort(G)^{\min}, \, \text{$\calF$ compatible}\}$, but their proof relies on another result \cite[Corollary~1]{louis}, which incorrectly assumes that minimal null supports correspond to minimal forts.  The Petersen graph provides a counterexample:  

\begin{ex} \label{ex:petersen}
Let $P$ be the Petersen graph. We have 
$\M(P) = \Z(P) = 5$~\cite{AIMMINRANK}, and so $\tau(\fort(P)^{\min}) = 5$ by Proposition~\ref{Z-tau}.  But we claim that any subcollection of $\fort(P)^{\min}$ that has transversal number equal to 5 (including $\fort(P)^{\min}$) is not compatible.

To prove our claim, we first note that the minimal forts of $P$ can be partitioned, up to isomorphism, into four distinct types, with five forts belonging to each type.  The  minimal forts of $P$ are as follows~\cite[Example~3.13]{cameron2023forts}:
   \begin{center}
   \begin{tabular}{c}
 $\{0, 1, 3, 8\}, \{0, 1, 9, 7\}, \{0, 2, 3, 5\}, \{0, 2, 4, 7\}, \{0, 8, 2, 9\}$, \\ 
 $\{0, 3, 6, 7\}, \{0, 8, 4, 6\}, \{0, 9, 5, 6\}, \{1, 2, 4, 9\},
 \{8, 1, 2, 5\}$,\\
 $\{1, 3, 4, 6\}, \{1, 3, 5, 9\}, \{8, 1, 4, 7\}, \{1, 5, 6, 7\}, \{9, 2, 3, 6\}$, \\
 $\{2, 4, 5, 6\},
 \{8, 2,6, 7\}, \{3, 4, 5, 7\},
\{8, 9, 3, 7\}, \{8, 9, 4, 5\}$.
\end{tabular}
\end{center} A representative of each type of fort is shown below: 
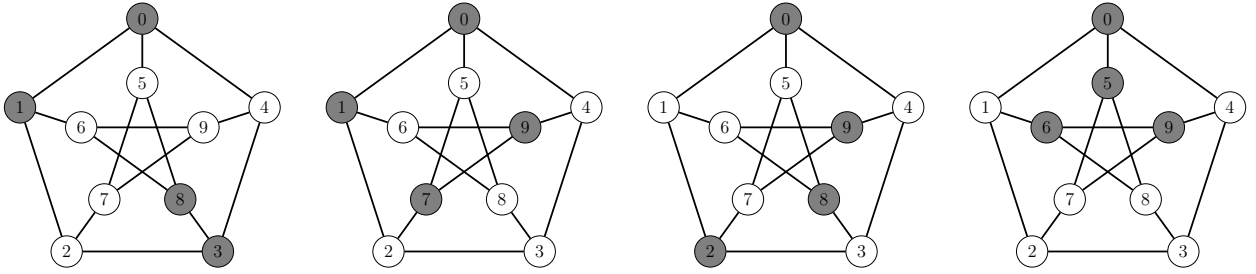
\begin{center}
\resizebox{\textwidth}{!}{\begin{tikzpicture}[scale=.75]
\node[circle,draw=black,fill=gray] (0) at (0,4) {$0$};
    \node[circle,draw=black,fill=gray] (1) at (-3.80,1.24) {$1$};
    \node[circle,draw=black,fill=white] (2) at (-2.35,-3.24) {$2$};
    \node[circle,draw=black,fill=gray] (3) at (2.35,-3.24) {$3$};
    \node[circle,draw=black,fill=white] (4) at (3.80,1.24) {$4$};
    \node[circle,draw=black,fill=white] (5) at (0,2) {$5$};
    \node[circle,draw=black,fill=white] (6) at (-1.90,0.62) {$6$};
    \node[circle,draw=black,fill=white] (7) at (-1.18,-1.62) {$7$};
    \node[circle,draw=black,fill=gray] (8) at (1.18,-1.62) {$8$};
    \node[circle,draw=black,fill=white] (9) at (1.90,0.62) {$9$};
    \foreach \x/\y in {0/1,1/2,2/3,3/4,4/0}
        \draw[black,=>latex',-,very thick] (\x) -- (\y);
    \foreach \x/\y in {0/5,1/6,2/7,3/8,4/9}
        \draw[black,=>latex',-,very thick] (\x) -- (\y);
    \foreach \x/\y in {5/7,7/9,9/6,6/8,8/5}
        \draw[black,=>latex',-,very thick] (\x) -- (\y);

\node[circle,draw=black,fill=gray] (10) at (0+10,4) {$0$};
    \node[circle,draw=black,fill=gray] (11) at (-3.80+10,1.24) {$1$};
    \node[circle,draw=black,fill=white] (12) at (-2.35+10,-3.24) {$2$};
    \node[circle,draw=black,fill=white] (13) at (2.35+10,-3.24) {$3$};
    \node[circle,draw=black,fill=white] (14) at (3.80+10,1.24) {$4$};
    \node[circle,draw=black,fill=white] (15) at (0+10,2) {$5$};
    \node[circle,draw=black,fill=white] (16) at (-1.90+10,0.62) {$6$};
    \node[circle,draw=black,fill=gray] (17) at (-1.18+10,-1.62) {$7$};
    \node[circle,draw=black,fill=white] (18) at (1.18+10,-1.62) {$8$};
    \node[circle,draw=black,fill=gray] (19) at (1.90+10,0.62) {$9$};

    \foreach \x/\y in {10/11,11/12,12/13,13/14,14/10}
        \draw[black,=>latex',-,very thick] (\x) -- (\y);
    \foreach \x/\y in {10/15,11/16,12/17,13/18,14/19}
        \draw[black,=>latex',-,very thick] (\x) -- (\y);
    \foreach \x/\y in {15/17,17/19,19/16,16/18,18/15}
        \draw[black,=>latex',-,very thick] (\x) -- (\y);

\node[circle,draw=black,fill=gray] (20) at (0+20,4) {$0$};
    \node[circle,draw=black,fill=white] (21) at (-3.80+20,1.24) {$1$};
    \node[circle,draw=black,fill=gray] (22) at (-2.35+20,-3.24) {$2$};
    \node[circle,draw=black,fill=white] (23) at (2.35+20,-3.24) {$3$};
    \node[circle,draw=black,fill=white] (24) at (3.80+20,1.24) {$4$};
    \node[circle,draw=black,fill=white] (25) at (0+20,2) {$5$};
    \node[circle,draw=black,fill=white] (26) at (-1.90+20,0.62) {$6$};
    \node[circle,draw=black,fill=white] (27) at (-1.18+20,-1.62) {$7$};
    \node[circle,draw=black,fill=gray] (28) at (1.18+20,-1.62) {$8$};
    \node[circle,draw=black,fill=gray] (29) at (1.90+20,0.62) {$9$};

    \foreach \x/\y in {20/21,21/22,22/23,23/24,24/20}
        \draw[black,=>latex',-,very thick] (\x) -- (\y);
    \foreach \x/\y in {20/25,21/26,22/27,23/28,24/29}
        \draw[black,=>latex',-,very thick] (\x) -- (\y);
    \foreach \x/\y in {25/27,27/29,29/26,26/28,28/25}
        \draw[black,=>latex',-,very thick] (\x) -- (\y);

\node[circle,draw=black,fill=gray] (30) at (0+30,4) {$0$};
    \node[circle,draw=black,fill=white] (31) at (-3.80+30,1.24) {$1$};
    \node[circle,draw=black,fill=white] (32) at (-2.35+30,-3.24) {$2$};
    \node[circle,draw=black,fill=white] (33) at (2.35+30,-3.24) {$3$};
    \node[circle,draw=black,fill=white] (34) at (3.80+30,1.24) {$4$};
    \node[circle,draw=black,fill=gray] (35) at (0+30,2) {$5$};
    \node[circle,draw=black,fill=gray] (36) at (-1.90+30,0.62) {$6$};
    \node[circle,draw=black,fill=white] (37) at (-1.18+30,-1.62) {$7$};
    \node[circle,draw=black,fill=white] (38) at (1.18+30,-1.62) {$8$};
    \node[circle,draw=black,fill=gray] (39) at (1.90+30,0.62) {$9$};

    \foreach \x/\y in {30/31,31/32,32/33,33/34,34/30}
        \draw[black,=>latex',-,very thick] (\x) -- (\y);
    \foreach \x/\y in {30/35,31/36,32/37,33/38,34/39}
        \draw[black,=>latex',-,very thick] (\x) -- (\y);
    \foreach \x/\y in {35/37,37/39,39/36,36/38,38/35}
        \draw[black,=>latex',-,very thick] (\x) -- (\y);
        
\end{tikzpicture}}
\end{center}

 For the sake of contradiction, assume that $\tau(\calF) = 5$ for some $\calF \subseteq \fort(P)^{\min}$.  Then $T = \{0,1,2,3,4\}$ is a minimum transversal of $\calF$.  The only forts that can satisfy Theorem ~\ref{correspondence} with respect to $T$ are the five of the last type, which must therefore belong to $\calF$.  Compatibility of $\{0,5,6,9\}$ and $\{2,6,7,8\}$ requires $\{0,2,8,9\}\in \calF$.  But then compatibility of $\{0,2,8,9\}$ and $\{0,5,6,9\}$ requires the existence of an element of $\fort(P)^{\min}$ contained in $\{0,2,5,6,8\}$, a contradiction.  Thus any subcollection of $\fort(P)^{\min}$ that has transversal number equal to 5 (including $\fort(P)^{\min}$) is not compatible.
\end{ex}

To rectify the issue of minimality, we can simply consider all forts.  Define  the \emph{fort transversal number} of $G$ as follows: \[ \Y(G) = \max \{\tau(\calF): \ \calF \subseteq \fort(G), \, \text{$\calF$ compatible}\}.\] Hicks et al.~\cite[proof of Theorem~1]{louis} showed that $\M(G) \leq \Y(G)$, but we can establish a more general result.

\begin{thm} \label{NYZ}
    For any graph $G$, $\N(G) \leq \Y(G) \leq \Z(G)$.
\end{thm}


\begin{proof}
    For the first inequality, first note that $\N(G) = \max \{\tau(\calFmin_A): \ A\in\csym(G)\}$ follows from Lemma \ref{tau(C_A)}.  For every $A\in\csym(G)$, $\calFmin_A$ is a compatible subcollection of $\fort(G)$ by Proposition \ref{csym forts} and Observation \ref{F_A^min-matroid}, so $\tau(\calFmin_A) = \tau(\calF_A) \leq \Y(G)$.
    By the comment following Proposition \ref{Z-tau}, $\Y(G) \leq \Z(G)$ is immediate.
\end{proof}

In general, $\M(G) \leq \N(G)$, and the next example shows these values can differ.  

\begin{ex}
    Minimal forts of the complete multipartite graph $K_{3,3,3}$ have either two vertices from one partite set or three vertices with one from each partite set.  The collection $\fort(K_{3,3,3})^{\min}$ can be seen to be compatible as the symmetric difference of any two minimal forts contains another minimal fort.  Moreover, $\tau(\fort(K_{3,3,3})^{\min}) = 7$.  Since $\M(K_{3,3,3}) = 6$ and $\Z(K_{3,3,3}) = 7$~\cite{FH14}, Theorem \ref{NYZ} implies $\M=6 < 7 = \Y=\Z$.  Berman~et~al.~\cite[Example~2.1]{berman2008mr} showed that $\N(K_{3,3,3}) \geq 7$, and hence $\N(K_{3,3,3}) = 7$.
\end{ex}

The next example reverses the strict and non-strict inequalities from the previous example.

\begin{ex}
Consider the corona of the 5-cycle with a single vertex, $G = C_5 \circ K_1$, shown below:
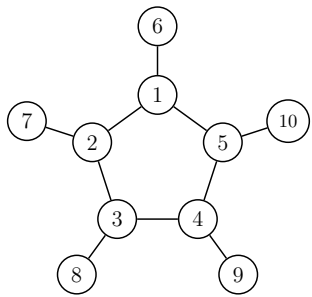
\begin{center}
 \resizebox{.25\textwidth}{!}{
        \begin{tikzpicture}[scale=1.3]
        \begin{scope}[every node/.style={circle,thick,draw}]
            \node[draw] (1) at (90:1) {$1$};
            \node[draw] (2) at (162:1) {$2$};
            \node[draw] (3) at (234:1) {$3$};
            \node[draw] (4) at (306:1) {$4$};
            \node[draw] (5) at (18:1) {$5$};
            \node[draw] (6) at (90:2) {$6$};
            \node[draw] (7) at (162:2) {$7$};
            \node[draw] (8) at (234:2) {$8$};
            \node[draw] (9) at (306:2) {$9$};
            \node[draw] (10) at (18:2) {\footnotesize $10$};
            \draw[thick] (5) -- (1) -- (2) -- (3) -- (4) -- (5); 
            \draw[thick] (1) -- (6);
            \draw[thick] (2) -- (7);
            \draw[thick] (3) -- (8);
            \draw[thick] (4) -- (9);
            \draw[thick] (5) -- (10);
        \end{scope}
        \end{tikzpicture}}
    \end{center}
    
It is known that $\M(G) = 2$ and $\Z(G) = 3$ \cite{BFH04, AIMMINRANK}.  We first claim that every fort of $G$ contains at least three pendant vertices.  Let $F$ be a fort of $G$.  If any vertex in the 5-cycle belongs to $F$, then its corresponding pendant vertex must also belong to $F$.  It is therefore impossible for $F$ to contain no pendant vertices. 
 It is similarly easy to see that $F$ cannot contain exactly one pendant vertex.  Suppose $F$ contains two pendant vertices, without loss of generality, $\{6,7\}$ or $\{6,8\}$.  If $1\in F$, then at least one more vertex from $\{4,5,9,10\}$ must belong to $F$, and so either $10$ or $9$
 must also belong to $F$.  Suppose $1\notin F$.  If $10 \notin F$, then $2\in F$.  If $\{6,7\} \subseteq F$, then either $8$ or $9$ to belong to $F$; if $\{6,8\} \subseteq F$, then $7$ must belong to $F$.  In all cases, $F$ must contain at least three pendant vertices.

 Next we claim that if $\calF$ is a collection of forts of $G$ and $T$ is a minimum transversal of $\calF$, then we may assume (without loss of generality) that $T$ consists only of pendant vertices.  Indeed, suppose $v\in T$ where $v$ is not a pendant vertex of $G$.  Then for each $F\in \calF$ such that $v \in F\cap T$, the pendant vertex $u$ adjacent to $v$ must also be contained in $F$.  So we may replace $v$ in $T$ by $u$ to get another minimum transversal of $\calF$. 

 Finally we show that $\tau(\calF) \leq 2$ for any compatible collection of forts $\calF$.  Suppose otherwise; that is, based on the above reasoning, suppose that for any two pendant vertices $\{x,y\}$, there exists a fort $F\in \calF$ such that $\{x,y\} \cap F = \emptyset$.  Then every set of exactly three pendant vertices must comprise the subset of pendant vertices of a fort in $\calF$.   In particular, $\calF$ must contain a fort $F$ whose three pendant vertices are $\{6,7,10\}$.  If $F$ contains $2$, then it must also contain $8$ or $9$, a contradiction.  So $2\notin F$, and similarly $5\notin F$.  It follows that $1\in F$ and that $F = \{1, 6, 7, 10\}$.  Similarly, $F' = \{2, 6, 7, 8\}$ must be a fort in $\calF$.  However, it is impossible for a fort to be contained in $F\cup F' \setminus \{6\} = \{1, 2, 7, 8, 10\}$.  Indeed, such a fort would need to be a subset of $\{2, 7, 8, 10\}$, which is impossible.

 Since there exists a compatible collection of forts with transversal number 2~\cite[Figure 1]{louis}, we see that $\Y(G) = 2$ and $\M = \N = \Y = 2 < 3 = \Z$.  
 \end{ex}

\begin{quest}
    Is there a graph $G$ with $\N(G) < \Y(G)$? 
\end{quest}

 By a result of Berman~et~al.~\cite[Corollary 2.6]{berman2008mr}, a graph with $\N(G) < \Y(G)$ must satisfy $\Z(G) > \delta(G)$; it must also have at least 8 vertices as $\M(G) = \Z(G)$ is known for all graphs up to order 7 \cite{DeLoss2010}, as well as several families of graphs \cite{AIMMINRANK}.  By Proposition \ref{theoneweremovedthenputback} such a graph would have the property that if $\calF \subseteq \fort(G)$ is a compatible collection of forts that achieves $\tau(\calF) = \Y(G)$, then there exists no $A \in \csym(G)$ such that $\calF_A = \calF$; indeed, it is straightforward to verify that $\calF^{\min}$ is a compatible collection, and if the matroid $\bm{M} = (V(G),\calFmin)$ were representable by a matrix $A\in \csym(G)$, then $\Y(G) = \tau(\calF) = \tau(\calF^{\min}) = n - \rank(\bm{M}) = n - \rank(A)$. 

\section{Disjoint Forts}\label{Section:DisjointForts}

One motivation for seeking to characterize $\M(G)$ and $\N(G)$ using transversal numbers of collections of forts is to gain a better understanding of how these two graph parameters differ. 
In this section, we prove that the requirement of symmetry can make a difference when looking for a matrix whose null vector supports are a given collection of forts, even when these forts are disjoint. 

\begin{prop}[{\cite[Proposition~6.3]{cameron2023forts}}]\label{disjoint}
    If $F_1, \dots, F_k$ are disjoint forts in $\fort(G)$, then there exists $A \in \csym(G)$ with $F_1,\dots,F_k \in \calF_A$.
\end{prop}

It is immediate that $\N(G)$ is at least the maximum number of disjoint forts that can be found in $G$, known as the {\em fort number} $\ft(G)$ \cite{cameron2023forts}.  However, we next show that Proposition~\ref{disjoint} is not true for symmetric matrices.  While it may still be true that $\M(G) \geq \ft(G)$~\cite[Question~6.1]{cameron2023forts}, this means a constructive proof for an arbitrary set of disjoint forts is not possible. 

We first need a definition:  For an undirected graph $G$, a {\em zero-sum flow} is an assignment of nonzero real
numbers to the edges such that, for each vertex, the sum of the values of all edges incident with
that vertex is zero.
According to Akbari~et~al.~\cite{zerosumtwo}, a connected bipartite graph has a zero-sum flow if and only if it is bridgeless, but that result is stated without proof.  Instead, a reference is provided (\cite{flow}), which only proves one direction.   As shown there,  it is possible for a graph with a bridge to have a zero-sum flow~\cite[Theorem 2 and Remark 2]{flow}.   So we provide a proof for arbitrary bipartite graphs for completeness. 

\begin{prop}\label{prop:questionmark}
If G is a bipartite graph, then G has a zero-sum flow if and only if $G$ is bridgeless.
\end{prop}

\begin{proof}
First note that $G$ has a zero-sum flow if and only if each of its connected components does, so we may assume that $G$ is connected.  

Let $G$ have a zero-sum flow and assume $G$ has a bridge for sake of contradiction.   Let $S_1$ and $S_2$ be the two partite sets of $G$.  Create a directed graph from $G$ by viewing each edge as starting in $S_1$ and ending in $S_2$.  Assigning the same edge values as in the zero-sum flow creates a flow in the more classical sense, which is known to be impossible in a graph with a bridge~\cite[Proposition 7.3.16]{west}. 

If $G$ is bridgeless, then every edge is part of a cycle.  Let $C_1,\dots, C_m$ be the simple cycles of $G$.  Since $G$ is bipartite, each $C_i$ is an even cycle. For each $C_i$, assign alternating values of $\pm \frac{1}{2^i}$ to its edges.  Then to each edge of $G$, assign the sum of the values from all cycles that contain that edge.  By construction, each edge is assigned a nonzero real number.  Since the sum of edge values at any vertex in each cycle is zero, the result is a zero-sum flow.
\end{proof}

\begin{thm} \label{bridgeless}
    Given a graph $G$, disjoint forts $F_1, F_2, \dotsc, F_m$ of $G$, and vectors $x_1, x_2, \dotsc, x_m$ with $\supp(x_i) = F_i$, there exists $A \in \msym(G)$ with $\{x_1,x_2,\dotsc,x_m\} \subseteq \Nul(A)$ if and only if the bipartite graph formed by the edges between $F_i$ and $F_j$ is bridgeless for each pair of forts. 
\end{thm}

\begin{proof}
  First, we can assume without loss of generality that the entries of each $x_i$ are in $\{0,1\}$.  Indeed, let $D_i$ be the diagonal matrix whose $j$th diagonal entry is $x_i(j)$ if $x_i(j)\neq 0$ and $1$ otherwise; then $y_i = D_i^{-1}x_i$ has entries in $\{0,1\}$.  Let $D = D_1 \cdots D_m$.  Then $A\in\msym(G)$ if and only if $\widehat A =DAD \in \msym(G)$.  Since the supports of the $x_i$ are disjoint, $Ax_i=0$ if and only if $\widehat A y_i = 0$.

Let $F =\bigcup_i F_i$.  For the forward direction, assume there exists $A \in \msym(G)$ with $\{x_1,x_2,\dotsc,x_m\} \subseteq \Nul(A)$.  Partition the matrix $A$ as 
\[
A=
\left[
\begin{matrix}
 C & B^T \\
B & E    
\end{matrix}
\right]
\]
where $C \in \msym(G[F])$ is further partitioned with respect to the $F_i$ and is of the form
\[
C=
\left[
\begin{matrix}
L_1 & M_{12}^T & \dotso & M_{1m}^T \\
M_{12} & L_2 & \dotso & M_{2m}^T \\
\vdots & \vdots & \ddots & \vdots \\
M_{1m} & M_{2m} & \dotso & L_m  
\end{matrix}
\right].
\]  

Since $Ax_i = 0$ for each $i$ and we have assumed the nonzero entries of $x_i$ are all equal to $1$,  we have $L_i{\bf 1}=0$,  $M_{ij}{\bf 1} = 0$, and $M_{ji}^T{\bf 1}=0$ for all $j \neq i$,  where ${\bf 1}$ is in each case an all-ones vector of the appropriate size.  Further, for each $i \neq j$, the matrix
\[ N_{ij}
 = \left[
\begin{matrix}
  0 & M_{ij}^T \\
M_{ij} & 0   
\end{matrix}
\right]
\]
has the zero/nonzero pattern of a bipartite graph $H_{ij}$.  The partite sets of $H_{ij}$ correspond to the (disjoint) supports of $x_i$ and $x_j$. Thus the entries of $M_{ij}$ yield a zero-sum flow on $H_{ij}$, which implies that  $H_{ij}$ has no bridges by Proposition \ref{prop:questionmark}.

Conversely, we will construct the symmetric matrix $A$, using the same notation as above, assuming that each bipartite graph $H_{ij}$ described by the matrix $N_{ij}$ is bridgeless.  By Proposition \ref{prop:questionmark}, $H$ has a zero-sum flow, so we may choose $M_{ij}$ so that $N_{ij}{\bf 1} = 0$.  For each $i$, let $L_i$ be the Laplacian matrix of $G[F_i]$.
Choose the entries of $B$ so that in each row of $B$, the nonzero entries corresponding to $F_i$ sum to zero (there are either no nonzero entries or at least two nonzero entries since $F_i$ is a fort),  and let $E$ be any matrix in $\msym(G([V(G)\setminus F])$.
\end{proof}

\begin{ex}
The graph below has disjoint forts \(F_1=\{1,4,6,7\}\) and \(F_2=\{2,3,5,8\}\) that have a bridge when the edges between them are considered as a bipartite graph, and so \(F_1\) and \(F_2\) cannot both be the supports of null vectors of a matrix \(A \in \msym(G)\) by Theorem \ref{bridgeless}.
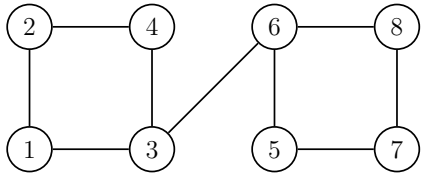
    \begin{center}
    \resizebox{.34\textwidth}{!}{
        \begin{tikzpicture}
        \begin{scope}[every node/.style={circle,thick,draw}]
            \node[draw] (1) at (0,0) {$1$};
            \node[draw] (2) at (0,2) {$2$};
            \node[draw] (3) at (2,0) {$3$};
            \node[draw] (4) at (2,2) {$4$};
            \node[draw] (5) at (4,0) {$5$};
            \node[draw] (6) at (4,2) {$6$};
            \node[draw] (7) at (6,0) {$7$};
            \node[draw] (8) at (6,2) {$8$};
            \draw[thick] (3) -- (1) -- (2) -- (4) -- (3) -- (6) -- (5) -- (7) -- (8) -- (6); 
        \end{scope}
        \end{tikzpicture}}
    \end{center}

    In particular a symmetric matrix of the form 
    
       \[
        \begin{bmatrix}
            a_{11} & a_{12} & a_{13} & 0 & 0 & 0 & 0 & 0 \\
            a_{12} & a_{22} & 0 & a_{24} & 0 & 0 & 0 & 0 \\
            a_{13} & 0 & a_{33} & a_{34} & 0 & a_{36} & 0 & 0 \\
            0 & a_{24} & a_{34} & a_{44} & 0 & 0 & 0 & 0 \\
            0 & 0 & 0 & 0 & a_{55} & a_{56} & a_{57} & 0 \\
            0 & 0 & a_{36} & 0 & a_{56} & a_{66} & 0 & a_{68} \\
            0 & 0 & 0 & 0 & a_{57} & 0 & a_{77} & a_{78} \\
            0 & 0 & 0 & 0 & 0 & a_{68} & a_{78} & a_{88} \\
        \end{bmatrix}   
    \]   
with two null vectors having supports $F_1 = \{1,4,6,7\}$ and $F_2 = \{2,3,5,8\}$, respectively, would force $a_{36}$ to be zero.  However, $F_1$ and $F_2$ are the supports of null vectors of the following combinatorially symmetric matrix: 

\[ C = 
        \begin{bmatrix}
            0 & -1 & 1 & 0 & 0 & 0 & 0 & 0 \\
            -1 & 0 & 0 & 1 & 0 & 0 & 0 & 0 \\
            1 & 0 & 0 & 1 & 0 & -2 & 0 & 0 \\
            0 & 1 & -1 & 0 & 0 & 0 & 0 & 0 \\
            0 & 0 & 0 & 0 & 0 & 1 & -1 & 0 \\
            0 & 0 & -2 & 0 & 1 & 0 & 0 & 1 \\
            0 & 0 & 0 & 0 & -1 & 0 & 0 & 1 \\
            0 & 0 & 0 & 0 & 0 & 1 & -1 & 0 \\
        \end{bmatrix}   
    \]
    
    On the other hand, forts \(F_3=\{1,4\}\), \(F_4 = \{5,8\}\) and \(F_5 = \{2,3,6,7\}\) are disjoint and no two have a bridge between them, so they can be simultaneous null vector supports for a symmetric matrix, which can be realized by replacing $c_{34}$ and $c_{78}$ by $-1$, and $c_{33}$ and $c_{66}$ by 2 in the matrix $C$ above.
\end{ex}

Disjoint forts have also played a recent role in the study of strong spectral properties of matrices.  A matrix $A\in\msym(G)$ has the Strong Arnold Property  (SAP) if the only symmetrix matrix $X$ satisfying $A \circ X = I \circ X = AX = 0$, where $\circ$ is the entrywise product and $I$ is the identity matrix, is the zero matrix.  Lin, Oblak,  and \v Smigoc~\cite[Lemma 2.3]{SSP} proved that if the vertices of a graph $G$ can be partitioned into three disjoint parts $R$, $W_1$, and $W_2$ that satisfy the following three properties (a \emph{barbell partition}), then there exists a matrix $A \in \msym(G)$ that does not have the SAP:
\begin{itemize}
\item $R$ can be empty but neither $W_1$ nor $W_2$ is empty;
\item there are no edges between $W_1$ and $W_2$;
\item for each $v \in R$ and $i \in \{1,2\}$, $|N_G(v) \cap W_i| \neq 1$. 
\end{itemize}
Allred~et~al.~noticed the connection of the third property to forts, defined two forts that have no edges between them to be \emph{separated}, and observed that $W_1$, $W_2$, and $R = V(G) \setminus (W_1 \cup W_2)$ is a barbell partition if and only if $W_1$ and $W_2$ are separated forts~\cite[Observation 2.5]{SSPop}.  We end by pointing out that Lin, Oblak,  and \v Smigoc's proof is similar to a special case of  Theorem~\ref{bridgeless}, taking $F_1 = W_1$, $F_2=W_2$, and $M_{12}=0$ in the construction of the matrix $A$.  Then $A \circ X = I \circ X = AX = 0$ where 
\[
X = \left[  
\begin{matrix}
0 & J & 0 \\ 
J^T & 0 & 0 \\
0 & 0 & 0
\end{matrix}
\right]
\]
where $J$ is a $|W_1|$-by-$|W_2|$ matrix of all ones.

\section*{Acknowledgments}

This project began as part of the ``Inverse Eigenvalue Problems for Graphs and Zero Forcing'' Research Community sponsored by the American Institute of Mathematics (AIM). We thank AIM for their support, and we thank the organizers and participants for contributing to this stimulating research experience.  Veronika Furst was supported in part by grant DMS-2331072 from the National Science Foundation.
\bibliographystyle{plainurl}
\bibliography{refs}

\end{document}